\documentclass[11pt,oneside,leqno]{amsart}
\allowdisplaybreaks
\usepackage[margin=2.5cm]{geometry}
\usepackage{multirow,tabularx}
\newcolumntype{Y}{>{\centering\arraybackslash}X}

\usepackage{pbox}
\usepackage{parskip}
\usepackage{eucal}
\usepackage{amscd}
\usepackage{amssymb}
\usepackage{amsmath}
\usepackage{amsfonts}
\usepackage{amsopn}
\usepackage{latexsym}
\usepackage{tikz}
\usepackage{comment}
\usepackage{hyperref}
\usepackage[nodisplayskipstretch]{setspace}
\usepackage{caption}
\usepackage{enumerate}
\usepackage[utf8]{inputenc}

\makeatletter
\@namedef{subjclassname@2010}{%
  \textup{2010} Mathematics Subject Classification}
\makeatother

\newtheorem{thm}{Theorem}[section]
\newtheorem{cor}[thm]{Corollary}
\newtheorem{lem}[thm]{Lemma}

\newtheorem{example}[thm]{Example}

\makeatletter
\newcommand{\BIG}{\bBigg@{2}}
\newcommand{\vast}{\bBigg@{3}}
\newcommand{\Vast}{\bBigg@{5}}
\makeatother

\numberwithin{equation}{section}

\begin{document}
\setlength{\arrayrulewidth}{0.1mm}



\title[Heron triangle and elliptic curves]{Heron triangles and a family of elliptic curves with rank zero}

\author[Vinodkumar Ghale]{Vinodkumar Ghale}
\address{Department of Mathematics\\ BITS-Pilani, Hyderabad campus\\
Hyderabad, INDIA}
\email{p20180465@hyderabad.bits-pilani.ac.in}

\author[Shamik Das]{Shamik Das}
\address{Department of Mathematics and Computation\\ IIT Guwahati\\
Guwahati, INDIA}
\email{shamikdas@iitg.ac.in}

\author[Debopam Chakraborty]{Debopam Chakraborty}
\address{Department of Mathematics\\ BITS-Pilani, Hyderabad campus\\
Hyderabad, INDIA}
\email{debopam@hyderabad.bits-pilani.ac.in}

\date{}


\maketitle

\section*{Abstract}

\noindent Given any positive integer $n$, it is well-known that there always exists a triangle with rational sides $a,b$ and $c$ such that the area of the triangle is $n$. For a given prime $p \not \equiv 1$ modulo $8$ such that $p^{2}+1=2q$ for a prime $q$, we look into the possibility of the existence of the triangles with rational sides with $p$ as the area and $\frac{1}{p}$ as $\tan \frac{\theta}{2}$ for one of the angles $\theta$. We also discuss the relation of such triangles with the solutions of certain Diophantine equations. 

\section{Introduction}

\noindent A triangle with rational side lengths and the rational area is called a Heron triangle. For a given integer $n$, the celebrated \textit{congruent number problem} is the problem of identifying a right angle triangle with $n$ as the area. The existence of $n$ as the area of a right angle triangle is equivalent to the rank of the congruent elliptic curve $E_{n}: y^{2} = x(x-n)(x+n)$ being positive. A glimpse of the extensive study of the congruent number problem can be found in \cite{Chahal}, \cite{Coates}, \cite{Johnstone}. In \cite{Goins}, the authors have proved that the existence of rational points of order greater than $2$ for a similar elliptic curve $E_{\Delta} : y^{2} = x(x-n \tau) (x+n \tau^{-1})$ is equivalent to the existence of a Heron triangle with the area as $n$. The elliptic curve associated with the Heron triangle appears in the works of various authors.\\
\noindent For a given integer $n$, the existence of infinitely many Heron triangle with area $n$ has been studied in \cite{Rusin}. The work of Buchholz and Ratbun in \cite{Buchholz} has proved the existence of infinitely many Heron triangles with two rational medians. Later in \cite{Buchholz2}, the authors have looked into the existence of Heron triangles with three rational medians. In \cite{Izadi}, the authors have shown the existence of elliptic curves associated with Heron triangle of rank at least three. Similar work has also been done to show the said rank to be at least four or five. In \cite{Halbeisen}, the authors have shown the existence of elliptic curves of rank at least two associated with Heron triangles.\\
\noindent The existence of prime pairs $(p,q)$ such that $p^{2}+1 = 2q$ is equivalent to the existence of a right-angled triangle with $p$ being one of the legs whose hypotenuse is also a prime number, namely $q$. It is conjectured that there are infinitely many such triangles. In this paper, we look into the Heron triangles with the area $p \equiv 5$ modulo $8$ such that $p^{2}+1 = 2q$ for some prime $q$. Moreover, we also assume the triangle to have one of the angles $\theta$ such that $\tau = \tan \frac{\theta}{2} = p^{-1}$. Using an explicit $2$-descent method, we calculate the rank of the elliptic curve $E_{\Delta_{ \tau,p}}: y^{2} = x(x-n \tau)(x+n \tau^{-1})$ associated to the above mentioned Heron triangle and then prove the following main result.

\begin{thm}\label{mainthm}
Let $p \equiv 5$ modulo $8$ be an odd prime such that $p^{2}+1 = 2q$ for a prime $q$. Then the elliptic curve $$E_{\Delta_{\tau,p}}: y^{2} = x(x-n \tau)(x+n \tau^{-1})$$ has rank always zero for $n =p$ and $\tau = \frac{1}{p}$. Hence, there exists no Heron triangle with area $p$ and one of the angles $\theta$ such that $\tan \frac{\theta}{2} = \frac{1}{p}$ whenever $p \equiv 5$ modulo $8$. Moreover, if $p \equiv 3,7$ modulo $8$, then the rank of $E_{\Delta_{\tau, p}}$ will be at most one.
\end{thm}
\noindent The works of Kramer and Luca in \cite{Kramer} showed a pathway towards the understanding of a Diophantine equation arising from different Heron triangles. We conclude this paper with a few observations towards the solvability of certain Diophantine equations associated to the above mentioned Heron triangles. In \cite{Yan}, X. Yan has shown that for fixed coprime positive integers $m$ and $n$ with different parity, the Diophantine equation $(m^{2}+n^{2})^{x} + (2mn)^{y} = (m+n)^{2z}$ has no solutions with $y \geq 2$. Here, we prove a one-one correspondence between a Heron triangle with area $p$ and $\tau = \frac{1}{p}$ and the Diophantine equation $(x^{2}+y^{2})^{2} + (2pxy)^{2} = z^{2}$. 

\section{\textbf{Heron Triangle and Elliptic Curve}}

\noindent We start this section with the relation between a Heron triangle and an elliptic curve. In \cite{Goins}, Goins and Maddox had shown that any triangle $\Delta_{n}$ with rational sides $a,b,c$, area $n \in \mathbb{Z}$ and one angle $\theta$ is associated with the following elliptic curve: $$E_{\Delta_{\tau, n}} : y^{2} = x(x - n \tau) (x + n \tau^{-1})$$ where $\tau$ denotes $\tan\frac{\theta}{2}$. Moreover, they have shown that the torsion group of $E_{\Delta_{\tau, n}}$ will be either $\mathbb{Z}/ 2 \mathbb{Z} \times \mathbb{Z}/ 2 \mathbb{Z}$ or $\mathbb{Z}/ 2\mathbb{Z} \times \mathbb{Z}/ 4\mathbb{Z}$, $\Delta$ being an isosceles triangle for the latter case. Using their result, for any pair $(p,q)$ of odd primes such that $p^{2}+1 = 2q$, one can identify the elliptic curve $E_{\Delta_{\tau, p}}$ with a Heron triangle $\Delta$ of area $p$ with one angle $\theta$ such that $\tau = \tan \frac{\theta}{2} = p^{-1}$ where $$E_{\Delta_{\tau, p}}: y^{2} = x(x-1)(x+p^{2}).$$  
By the method of 2-descent (see \cite{Silverman}), there exists an injective homomorphism $b$ such that
$$ b: E_{\Delta_{\tau, p}}(\mathbb{Q})/2E_{\Delta_{\tau, p}}({\mathbb{Q}}) \longrightarrow \mathbb{Q}(S,2) \times \mathbb{Q}(S,2)$$
where\\ $\mathbb{Q}(S,2) = \{b \in \mathbb{Q}^{*}/(\mathbb{Q}^{*})^{2}\text{ : } ord_{v}(b) \equiv 0 \text{ (mod 2) for } v \neq 2,p,q\} = \{ \pm 1,\pm 2,\pm p, \pm q, \pm 2p, \pm 2q, \pm pq, \pm 2pq\}$. The map is defined by 
\begin{align*}
b(x,y) = \begin{cases}
(x, x -1)  & \text{if} \text{  } x\neq 0,1. \\
(-1,-1) & \text{if} \text{  } x = 0.\\
(1,2q) & \text{if}\text{  } x =1. \\
(1,1) &\text{if} \text{  } x = \infty \text{  } i.e.\text{  } P = \mathcal{O}.
\end{cases}
\end{align*}
Moreover, if $(b_1, b_2) \in \mathbb{Q}(S,2) \times \mathbb{Q}(S,2)$ is a pair that is not in the image of one of the three points $O , (0,0), (1,0)$, then $(b_1,b_2)$ is the image of a point $ P = (x,y) \in E_{\Delta_{\tau,p}} (\mathbb{Q})/2E_{\Delta_{\tau,p}}(\mathbb{Q})$ if and only if the equations 
\begin{eqnarray} b_1z_1^2 - b_2z_2^2 = 1 \label{eq22} \\ 
 b_1z_1^2 - b_1b_2z_3^2 = -p^2 \label{eq23}
\end{eqnarray}
have a solution $(z_1,z_2,z_3) \in K^* \times K^* \times K$. If $(z_{1},z_{2},z_{3})$ is  solution to the equations (\ref{eq22}) and (\ref{eq23}), then the pre-image $P = (x,y) \in E_{\Delta_{\tau,p}}$ of a point $(b_{1},b_{2})$ under the map $b$ is given by $x = b_{1}z_{1}^{2}, \quad y = b_{1}b_{2}z_{1}z_{2}z_{3}$.\\ 
\noindent Reducing $E_{\Delta_{\tau, p}}$ modulo $3$ for $p \neq 3$ while simultaneously noticing the fact $q \not \equiv 0$ (mod $3$), one can immediately observe that $|\Tilde{E}_{\Delta_{\tau, p}}(\mathbb{F}_{3})| = 4$. Since $E[2] \subset E(\mathbb{Q})_{tors}$ and $E(\mathbb{Q})_{tors}$ injects into $\Tilde{E}(\mathbb{F}_{3})$ for any elliptic curve $E$, one can see that $E_{\Delta_{\tau,p}}(\mathbb{Q})_{tors} = E_{\Delta_{\tau,p}}[2] \text{ whenever p} \neq 3$. A similar approach shows that $E_{\Delta_{\tau,3}}(\mathbb{Q})_{tors} \cong \mathbb{Z}/ 2 \mathbb{Z} \times \mathbb{Z}/ 2 \mathbb{Z}$ too. Hence we can start with the following observation: $$E_{\Delta_{\tau, p}}(\mathbb{Q})_{tors} \cong \mathbb{Z}/ 2 \mathbb{Z} \times \mathbb{Z}/ 2 \mathbb{Z}.$$ We first start with the following result regarding the modifications of the equations (\ref{eq22}) and (\ref{eq23}). 
\begin{lem}
Let $(b_1, b_2) \in \mathbb{Q}(S,2) \times \mathbb{Q}(S,2)$ is a pair that is not in the image of points of $E_{\Delta_{\tau,p}}(\mathbb{Q})_{tors}$. Then $(b_1,b_2)$ is the image of a point
$ P = (x,y) \in E_{\Delta_{\tau, p}}(\mathbb{Q})/2E_{\Delta_{\tau, p}}(\mathbb{Q})$ if and only if the equations \begin{equation}\label{maineq1}
     b_1r_1^2 - b_2r_2^2 = s^{2},
\end{equation}
\begin{equation}\label{maineq2}
 b_1r_1^2 - b_1b_2r_3^2 = -p^2s^2,
\end{equation}
have a solution $(r_1,r_2,r_3,s) \in \mathbb{Z}^* \times \mathbb{Z}^* \times \mathbb{Z}^{\ast} \times \mathbb{Z}^* $ with $(r_{1},s)=(r_{2},s)=(r_{3},s)=1$.
\end{lem}
\noindent \textbf{Proof:} Assuming $z_i = \dfrac{r_i}{s_i}$ for $i = 1,2,3$ where $ r_i, s_i \in \mathbb{Z}$ and $ (r_i, s_i) =1$ for $i = 1,2,3$, the result follows immediately if one assumes $(b_1,s_1) = (b_2,s_2) = (b_1b_2,s_3) = 1$ in addition. In case the additional assumption of $(b_1,s_1) = (b_2,s_2) = (b_1b_2,s_3) = 1$ is not true, a simple calculation for each of the three cases here will yield into contradiction. We show the proof for one of those three cases, the other two follow a similar approach towards their solution. If $(b_1, s_1) \neq 1 $ then there exists a prime number $t_1$ such that $t_1|b_1$ and $t_1|s_1$. Now, from the equation (\ref{eq22}), we have  $r_1^2s_2^2 \equiv 0 \ (mod \ t_1)$ as $b_1$ is square-free. This in turn implies that $s_2 \equiv 0 \ (mod  \ t_1)$ as $(r_1,s_1) = 1 $ and $t_1 | s_1$. If $ t_1 \nmid b_2$, then $b_{1}r_{1}^{2}s_{2}^{2}$ is divisible by an odd power of $t_{1}$ whereas $b_2r_2^2s_1^2+s_1^2s_2^2$ is divisible by an even power of $t_{1}$, a contradiction because $b_{1}r_{1}^{2}s_{2}^{2} = b_2r_2^2s_1^2+s_1^2s_2^2$. Similarly if $t_{1} | b_{2}$, then an odd power of $t_{1}$ divides $b_1r_1^2s_2^2 - b_2r_2^2s_1^2$ whereas an even power of $t_{1}$ divides $s_1^2s_2^2$, a contradiction again.  $\hfill \square$

\noindent In total there are $256$ different possibilities of $(b_{1}, b_{2}) \in \mathbb{Q}(S,2) \times \mathbb{Q}(S,2)$. The next result shows that it is sufficient to focus only on $16$ of those pairs. Before we state our next result, for the sake of brevity we will assume $A$ is the image of the set $E_{\Delta_{\tau,p}}(\mathbb{Q})_{tors}$ under the map $b$ that is, $$A = \{(-1,-1), (1,2q), (1,1), (-1, -2q)\}.$$
\begin{lem}\label{lemnewest}
Suppose $(b_{1}, b_{2}) \in \mathbb{Q}(S,2) \times \mathbb{Q}(S,2)$ such that $(b_{1}, b_{2})$ is image point of the group $E_{\Delta_{\tau,p}}(\mathbb{Q})$ modulo $2E_{\Delta_{\tau,p}}(\mathbb{Q})$. Then:\\
$[a]$ $b_{1}b_{2} > 0$ always.\\
$[b]$ $b_{1}$ is odd always.\\
$[c]$ If $(b_{1},b_{2}) \in Im(b)$ modulo $A$, then  $b_{i} \in \{1,p,q,pq\}$ for $i=1,2$.
\end{lem}
\noindent \textbf{Proof:}  $b_{1}b_{2}<0$ together with equations (\ref{maineq1}) and (\ref{maineq2}) implies either $p^{2} < 0$ or $s^{2} < 0$, a contradiction in either case. This immediately proves $[a]$.\\
We now show that $b_{1}$ always has to be odd. Because $p$ is an odd prime, from the equation (\ref{maineq2}) one can get  $p^2s^2 = b_1b_2r_3^2 - b_1r_1^2 \equiv 0 \ (mod \ 2) \implies s \equiv 0$ $(mod \ 2)$ if $b_{1}$ is even.
Then from the equation (\ref{maineq1}) we get, $b_1r_1^2 - b_2r_2^2 = s^2 \equiv 0 \ (mod \ 4)$, which in turn implies $b_2 \equiv b_1 \equiv 0 \ (mod \ 2) $. This is because $ ( r_1,s) = 1 = (r_2,s)$ and $s \equiv 0 \ (mod \ 2 ) $. So, $ r_1 \equiv r_2 \equiv 1 \ (mod \ 2 )$. Now from the equation (\ref{maineq2}) we get, 
\begin{equation*}
    b_1r_1^2 = b_1b_2r_3^2 - p^2s^2 \equiv 0 \ (mod \ 4) \\ \implies r_1 \equiv 0 \ (mod \ 2)
\end{equation*} as $b_1$ is square-free. This contradicts the fact that $(r_1, s_1) = 1 $. Hence, $b_1 \not \equiv 0 \ (mod \ 2)$ and the assertion of part $[b]$ holds.\\
\noindent Now for part $[c]$, using the fact that $b$ is a homomorphism, we can now immediately say that any pair $(b_{1},b_{2}) \in Im(b)$ if and only if $(b_{1},b_{2}) \times (a_{1},a_{2}) = (a_{1}b_{1}, a_{2}b_{2}) \in Im (b)$ where $(a_{1},a_{2}) \in A$. Hence without loss of generality, one can only focus on examining the possibility of $(b_{1},b_{2}) \in Im(b)$ where both $b_{1}$ and $b_{2}$ belong to $Im(b)/A = \{1,p,q,pq\}$. This proves the statement of part $[c]$. $\hfill \square$ 

\noindent We now systematically look into the rest of the sixteen possibilities for $(b_{1},b_{2})$ as potential image points under the map $b$. The next result shows there will be only one point to be concerned with.
\begin{lem}
Let $(b_{1},b_{2})$ be a pair such that $b_{i} \in \{1,p,q,pq\}$ for $i=1,2$. If $(b_{1},b_{2}) \neq (1,q)$ then $(b_{1},b_{2}) \not \in Im(b)$ whenever $p \equiv 5$ (mod $8$).
\end{lem}
\noindent \textbf{Proof:} Suppose gcd $(b_{1},b_{2}) \neq 1$ for some $(b_{1},b_{2}) \in Im(b)$. Then gcd $(b_{1},b_{2})$ is $p,q$ or $pq$. If $p$ divides gcd $(b_{1},b_{2})$ then $s \equiv 0$ (mod $p$) from equation (\ref{maineq1}). Also from equation (\ref{maineq2}), one can observe that $ b_1b_2r_3^2 - p^2s^2 = b_1r_1^2 \equiv 0$ modulo $p^{2}$. Hence, $r_{1} \equiv 0$ modulo $p$ as $b_{1}$ is square-free. This leads us to a contradiction as then $p$ divides gcd $(r_1,s) = 1$. So $(b_{1},b_{2}) \not \in Im(b)$ if gcd $(b_{1},b_{2}) = p$ or $pq$. A very similar argument proves that $q$ can also not divide gcd $(b_{1},b_{2})$ if $(b_{1},b_{2}) \in Im(b)$. So $(b_{1},b_{2}) \in Im(b)$ implies that gcd $(b_{1},b_{2}) = 1$. Now we are left with eight possible pairs of $(b_{1},b_{2})$, i.e. $(1,p), (1,q), (1,pq), (p,1), (p,q), (q,1), (q,p)$ and $(pq,1)$ (ignoring $(b_{1},b_{2}) = (1,1)$ which is the image of $\mathcal{O}$ under the map $b$).\\
\noindent If $(b_{1},b_{2}) = (1,pq) \in Im(b)$ then from equation (\ref{maineq2}), one gets that $pqr_{3}^{2}-p^{2}s^{2} = r_{1}^{2} \equiv 0 \ (mod \ p)$. But then from equation (\ref{maineq1}), $p$ also divides $r_{1}^{2} - pqr_{2}^{2} = s^{2}$, a contradiction as gcd $(r_{1},s) = 1$. So $(1,pq) \not \in Im(b)$. A very similar argument shows that $(pq,1) \not \in Im(b)$ either.\\
\noindent $(p,q) \not \in Im(b)$ as $(p,q) \in Im(b)$ implies $pr_{1}^{2} - pqr_{3}^{2} = -p^{2}s^{2}$ from equation (\ref{maineq2}). This in turn implies that $2r_{1}^{2} \equiv 2qr_{3}^{2} \equiv r_{3}^{2}$ (mod $p$), hence either $\Big(\frac{2}{p}\Big) = 1$ or $r_{1} \equiv r_{3} \equiv 0$ (mod $p$), a contradiction either way as $p \equiv 5$ (mod $8$) and $(r_{1},s) = 1$. If $(q,p) \in Im(b)$, then from equation (\ref{maineq2}), we have $qr_1^2 - pqr_3^2 = -p^2s^2   \implies r_1^2 \equiv 0 \ (mod \ p)$. This leads to a contradiction because from equation (\ref{maineq1}), $s^{2} = qr_{1}^{2} - pr_{2}^{2} \implies s \equiv 0$ (mod $p$), hence $p$ divides both $r_{1}$ and $s$. So $(q,p) \not \in Im(b)$. Now we are left with only four possible image points for the homomorphism $b$. Those point are $(p,1)$ $(1,p), (q,1)$ and $(1,q)$.\\
\noindent $(q,1) \not \in Im(b)$ because otherwise $s \equiv 0$ (mod $q$) from equation (\ref{maineq2}) and then $r_{2} \equiv 0$ (mod $q$) from equation (\ref{maineq1}). This leads to a contradiction as gcd $(r_{2},s) = 1$. Now $(1,p) \in Im(b)$ implies that $r_1 \equiv 0 \ (mod \ p)$ from equation (\ref{maineq2}). Then from equation (\ref{maineq1}) one can observe that $r_1^2 - pr_2^2 = s_2^2 \equiv 0 \ (mod \ p)$ which implies $s \equiv 0 \ (mod \ p)$, a contradiction, as $ (r_1, s) = 1$. Similarly $(p,1) \in Im(b)$ implies $pr_3^2 - r_2^2 = 2qs^2 \implies p r_{3}^{2} \equiv r_{2}^{2}$ (mod $q$) from equations (\ref{maineq1}) and (\ref{maineq2}), a contradiction if $r_{2} \not \equiv 0$ (mod $q$) as $\Big(\frac{p}{q}\Big) \neq 1$ from Gauss' quadratic reciprocity law and the fact $p^{2}+1 = 2q$. But if $r_{2} \equiv r_{3} \equiv 0$ (mod $q$), then again from the above equation, $2qs^{2} = p r_{3}^{2} -r_{2}^{2} \equiv 0$ (mod $q^{2}$) which implies $s \equiv 0$ (mod $q$), a contradiction as gcd $(r_{2},s) = 1$. Hence $(1,p), (p,1) \not \in Im(b)$ and the result follows. $\hfill \square$ \\
\noindent We have excluded the possibility of every non-trivial pair $(b_{1},b_{2})$ being in the image of the homomorphism $b$, apart from the point $(1,q)$. now we prove the following lemma which essentially proves the Theorem \ref{mainthm}.
\begin{lem}
$(1,q) \not \in Im(b)$ if $p \equiv 5$ (mod $8$).
\end{lem}
\noindent \textbf{Proof:} From the equations (\ref{maineq1}) and (\ref{maineq2}) we know that if $(1,q) \in Im (b)$, then 
\begin{eqnarray}
    r_1^2 - qr_2^2 = s^2, \label{newmaineq1} \\
    r_1^2 - qr_3^2 = -p^2s^2. \label{newmaineq2}
\end{eqnarray} for some $(r_1,r_2,r_3,s) \in  \mathbb{Z^*} \times \mathbb{Z^*} \times \mathbb{Z}^{*} \times \mathbb{Z^*} $. This in turn implies $r_3^2 - r_2^2 = qs^2$. From equation (\ref{newmaineq2}), we get that either $\Big(\frac{q}{p}\Big) = 1$ or $r_{1} \equiv r_{3} \equiv 0$ (mod $p$). But $p \equiv 5$ (mod $8$) and $2q = p^{2}+1$ implies that $\Big(\frac{q}{p}\Big) \neq 1$ and hence $r_{1} \equiv r_{3} \equiv 0$ (mod $p$). But from equation (\ref{newmaineq1}), $r_{1} \equiv 0$ (mod $p$) implies either $-q \equiv 0$ (mod $p$) or $r_{2} \equiv s \equiv 0$ (mod $p$), a contradiction in either way as $\Big(\frac{-q}{p}\Big) \neq 1$ and gcd $(r_{2},s) = 1$. Hence $(1,q) \not \in Im(b)$.   $\hfill \square$ \\

\noindent Now for $p \equiv 5$ (mod $8$), the result of the Theorem \ref{mainthm} follows directly from the observation that $(b_{1},b_{2}) \in Im(b)$ only when they are the images of one of the four torsion points. Hence the rank of the elliptic curve $E_{\Delta_{\tau,p}}$ is zero when $p \equiv 5$ modulo $8$. We enlist a table of few such curves and their corresponding rank below in support of our claim.\\
(Corresponding elliptic curve is $E: y^2 = x^3 + (p^2 -1)x^2 - p^2 x$ and the rank is verified @ http://magma.maths.usyd.edu.au/calc/ ).\\
\noindent For $p \equiv 3,7$ (mod $8$), we first observe that the same proof for $p \equiv 5$ (mod $8$) works, except for the image points $(b_{1},b_{2}) = (1,q), (p,q)$ and $(p,1)$ where we have used the fact $p \equiv 5$ (mod $8$). In the following lemma we prove that two of those three points can not appear as image points even for the cases $p \equiv 3,7$ (mod $8$).\\
\begin{lem}
If $p \equiv 3$ (mod $4$), then $(p,q), (p,1) \not \in Im(b)$.  
\end{lem}
\noindent \textbf{Proof:} First let us suppose that $p \equiv 3$ (mod $8$). Then $(p,q) \not \in Im(b)$ as $(p,q) \in Im(b)$ implies $pr_{1}^{2} - pqr_{3}^{2} = -p^{2}s^{2}$ from equation (\ref{maineq2}). This in turn implies that $2r_{1}^{2} \equiv 2qr_{3}^{2} \equiv r_{3}^{2}$ (mod $p$), hence either $\Big(\frac{2}{p}\Big) = 1$ or $r_{1} \equiv r_{3} \equiv 0$ (mod $p$), a contradiction either as $\Big(\frac{2}{p}\Big) = -1$ whenever $p \equiv 3$ (mod $8$) and $r_{1} \equiv r_{3} \equiv 0$ implies that $-p^{2}s^{2} \equiv 0$ (mod $p^{3}$) from equation (\ref{maineq2}) which in turn implies $s \equiv 0$ (mod $p$), a contradiction as gcd $(r_{1},s) = 1$. Similarly, when $p \equiv 7$ (mod $8$), because from equation (\ref{maineq1}), one gets $r_{2}^{2} \equiv -2s^{2}$ (mod $p$) which implies $r_{2} \equiv s \equiv 0$ (mod $p$) or $\Big(\frac{-2}{p}\Big) = 1$, contradiction either way as gcd $(r_{2},s)=1$ and $\Big(\frac{-2}{p}\Big) = -1$ when $p \equiv 7$ (mod $8$). Hence $(p,q) \not \in Im(b)$ if $p \equiv 3$ (mod $4$).\\
\noindent Now $(p,1) \in Im(b)$ implies that $pr_{1}^{2} - r_{2}^{2} = s^{2}$ from equation (\ref{maineq1}). This in turn implies $r_{2}^{2} \equiv -s^{2}$ (mod $p$). Hence either $\Big(\frac{-1}{p}\Big) = 1$ or $r_{2} \equiv s \equiv 0$ (mod $p$), again contradiction either way as $p \equiv 3$ (mod $4$) and gcd $(r_{2},s) = 1$. So $(p,1) \not \in Im(b)$ when $p \equiv 3$ (mod $4$). $\hfill \square$ \\
\noindent Now we can conclude the Theorem \ref{mainthm} by observing that for the cases $p \equiv 3,7$ (mod $8$), there exists only one possible image point (modulo the image of the torsion group), that is $(1,q)$. Hence, when $p \equiv 3,7$ (mod $8$), the rank of $E_{\Delta_{\tau,p}}$ is at most $1$. We enlist a few examples below of such curves with rank one.
\begin{center}
\label{t}
\begin{tabular}{ |c|c|c|c|c|c|c|c| }
\hline

 p & p mod 8 &  q & Rank & p & p mod 8 & q & Rank \\ 
 \hline
 3 & 3 & 5 & 1 &  61 & 5 & 1861 & 0\\ 
  5 &	5 & 13  & 0 & 71 & 7 & 2521 & 1 \\
  11 & 3 & 61 & 1 & 79 & 7 & 3121 & 1 \\
29	& 5 &  421 & 0 & 739 & 3 & 273061 & 1\\ 
\hline
\end{tabular}
\end{center}


\section{\textbf{Diophantine Equation and Heron Triangle}}

\noindent Given a Heron triangle with area $p$ and one of the angles $\theta$ such that $\tau = \tan \frac{\theta}{2} = \frac{1}{p}$, in this section we show that the Diophantine equation $(x^{2}+y^{2})^{2} + (2pxy)^{2} = z^{2}$ with gcd $(x,y) = 1$ is solvable if and only if there exists a Heron triangle $\Delta_{\tau,p}$ of the form mentioned above. We then conclude that there is no solution to the aforementioned Diophantine equation when $p \equiv 5$ modulo $8$ with the help of the Theorem \ref{mainthm} from the previous section. The main result of this section is as follows:
\begin{thm}\label{mainthm1}
For every fixed odd prime $p$, there is a one-one correspondence between the solvability of the Diophantine equation $(x^{2}+y^{2})^{2} + (2pxy)^{2} = z^{2}$ with gcd $(x,y) = 1$ and the existence of a Heron triangle with area $p$ and one of the angles $\theta$ such that $\tan \frac{\theta}{2} = \frac{1}{p}$. Moreover, if $p \equiv 5$ modulo $8$, then $(x^{2}+y^{2})^{2} + (2pxy)^{2} = z^{2}$ has no solution if $p^{2}+1 = 2q$ for some prime $q$.
\end{thm}
\noindent \textbf{Proof of Theorem \ref{mainthm1}:}
\noindent We first notice that both $\sin \theta$ and $\cos \theta$ will be rational. From the laws of sines and cosines for a triangle, we get the following;
\begin{equation}\label{eq18}
    \cos \theta = \frac{a^{2}+b^{2} - c^{2}}{2ab},  \text{ }  \sin \theta  = \frac{2p}{ab}, \text{ } \tan \frac{\theta}{2} = \frac{4p}{(a+b)^{2} - c^{2}}.
\end{equation}
\noindent From the equation (\ref{eq18}) and the fact that we have assumed $\tan \frac{\theta}{2} = \frac{1}{p}$, we arrived at the following conclusion; 
\begin{equation}\label{eq19}
    (a+b)^{2} - c^{2} = 4p^{2} \implies (a+b)^{2} = c^{2} + (2p)^{2}.
\end{equation}
Also from the equation (\ref{eq18}) and the fact that $\sin^{2} \theta + \cos^{2} \theta = 1$, we get the following;
\begin{equation}\label{eq20}
    \Big(\frac{a^{2}+b^{2} - c^{2}}{2ab}\Big)^{2} + \Big(\frac{2p}{ab}\Big)^{2} = 1 \implies ab = 1+ p^{2}. 
\end{equation}
\noindent Equations (\ref{eq19}) and (\ref{eq20}) together implies that $(a-b)^{2} = (a+b)^{2} - 4ab = c^{2}-4$. Assuming $a-b = u = \frac{u_{1}}{u_{2}}$ and $c = w = \frac{w_{1}}{w_{2}}$, where both the representations are in their lowest forms, we observe that $\frac{w_{1}^{2} - 4w_{2}^{2}}{w_{2}^{2}} = \frac{u_{1}^{2}}{u_{2}^{2}}$. Because gcd $(w_{1}^{2} - 4w_{2}^{2}, w_{2}^{2}) = 1$, we can write $u_{2}^{2} = w_{2}^{2}$ and $u_{1}^{2} = w_{1}^{2} - 4w_{2}^{2}$ which then implies $w_{1}^{2} = u_{1}^{2} + 4w_{2}^{2}$. From the fact that gcd $(w_{1}, w_{2}) = 1$, it can immediately be seen $w_{1}, 2w_{2}$ and $u_{1}$ are all pairwise coprime and they form a \textit{Pythagorean primitive triplet}. Hence, there exist integers $m,n$ with $m > n$ such that gcd $(m,n) = 1$ and $u_{1} = m^{2} - n^{2}, \text{ } 2w_{2} = 2mn \text{ and } w_{1} = m^{2}+n^{2}.$
This implies $c = \frac{w_{1}}{w_{2}} = \frac{m^{2}+n^{2}}{mn} = \frac{m}{n}+\frac{n}{m}$. Similarly from the fact that $(a-b)^{2} = \frac{u_{1}^{2}}{w_{2}^{2}}$, one can observe that $a-b = \frac{m}{n} -  \frac{n}{m}$, under the assumption $a \geq b$ which implies that $b = a - (\frac{m}{n} - \frac{n}{m})$. If $b \geq a$, one can change $a-b$ suitably. Now from the equation (\ref{eq19}), we know that $(a+b)^{2} = 4p^{2} + c^{2} = 4p^{2} + (\frac{m}{n} + \frac{n}{m})^{2}$, where replacing $b$ as in the previous line, we get the following;
\begin{equation}\label{eq21}
     a = \frac{\frac{m}{n} - \frac{n}{m} \pm \sqrt{(\frac{m}{n} + \frac{n}{m})^{2}+4p^{2}}}{2} = \frac{m^{2} - n^{2} \pm \sqrt{(m^{2}+n^{2})^{2} + (2pmn)^{2}}}{2mn}.
\end{equation}
This immediately gives us a solution to the Diophantine equation $z^{2} = (x^{2}+y^{2})^{2} + (2pxy)^{2}$ with gcd $(x,y) = 1$ in the form of 
\begin{equation}\label{eq00}
z = 2mna - (m^{2} - n^{2}), \text{ } x = m \text{ and } y = n.    
\end{equation}
\noindent For the converse implication, we start with the assumption that there exists positive integers $x,y$ and $z$ such that gcd $(x,y) = 1$ and $z^{2} = (x^{2}+y^{2})^{2} + (2pxy)^{2}$. Define 
\begin{equation}\label{eq01}
a = \frac{\frac{x}{y} - \frac{y}{x} + \sqrt{(\frac{x}{y}+\frac{y}{x})^{2}+4p^{2}}}{2} = \frac{\frac{x}{y} - \frac{y}{x} + \frac{z}{xy}}{2}, \text{ }b = a - (\frac{x}{y} - \frac{y}{x}) \text { and } c = \frac{x}{y} + \frac{y}{x}.    
\end{equation}
Then $a,b,c$ are all positive integers. Also, $ab = 1+p^{2} \text{ and } a-b = \frac{y}{x} - \frac{x}{y}.$ Now suppose we have a triangle with sides $a,b$ and $c$ and the angle $\theta $ between the sides of the length $a$ and $b$. Then by the law of cosines, $$\cos \theta = \frac{a^{2}+b^{2} - c^{2}}{2ab} = \frac{(a-b)^{2}-c^{2}}{2ab}+1 = \frac{p^{2}-1}{p^{2}+1} \in \mathbb{Q}.$$
Similarly, from the formula of $\sin \theta  = \sqrt{1 - \cos ^{2} \theta}$ we get, $\sin \theta  = \frac{2p}{p^{2}+1} = \frac{2p}{ab}$. This implies that the triangle have area $\frac{1}{2}ab \sin \theta = p$. Substituting the values of $a,b$ and $c$ in terms of $x,y$ and $z$ and observing the fact that $z^{2} = (x^{2}+y^{2})^{2} + (2pxy)^{2}$, one can easily observe that $\tan \frac{\theta}{2} = \frac{4p}{(a+b)^{2} - c^{2}} = \frac{1}{p}$ too. This concludes the proof of the statement. \hfill $\square$

\noindent The following result looks into the solution of the Diophantine equation $x^{2}+y^{2}+x^{2}y^{2} = p^{2}$ for some odd prime number $p$. The result follows directly from Theorem \ref{mainthm1}.
\begin{cor}\label{cor1}
For an odd prime $p$, there exists a Heron triangle with area $p$ and one rational angle $\theta$ such that $\tan \frac{\theta}{2} = \frac{1}{p}$ whenever the Diophantine equation $p^{2} = x^{2}+y^{2}+x^{2}y^{2}$ is solvable. 
\end{cor}
\noindent \textbf{Proof:} By Pythagoras' primitive element theorem, any solution $(m,n,l)$ of $(x^{2}+y^{2})^{2} + (2pxy)^{2} = z^{2}$ with gcd $(m,n) = 1$ implies that there exists co-prime natural numbers $r$ and $s$ such that $m^{2}+n^{2} = r^{2} - s^{2}$, $2pmn = 2prs$ and $l = r^{2}+s^{2}$. Hence $r = p$ and $s = mn$ will be a possible solution if $r^{2} - s^{2} = p^{2} - m^{2}n^{2} = m^{2}+n^{2}$. This concludes the proof of the statement.
\begin{example}
For $p=3$, one can notice that $x=2, y=1$ satisfies $x^{2}+y^{2}+x^{2}y^{2} = 9$. Using the explicit formula for $a,b,c$ from the Theorem \ref{mainthm1}, one can easily obtain $a = 4$, $b = \frac{5}{2}$ and $c = \frac{5}{2}$ to be the sides of a triangle with area $3$ and one angle $\theta$ such that $\sin \theta, \cos \theta \in \mathbb{Q}$ and $\tan \frac{\theta}{2} = \frac{1}{3}$. \hfill $\square$
\end{example}

\section*{Acknowledgement} 
\noindent The authors would like to thank Prof. Anupam Saikia for providing valuable suggestions throughout this work. The first author would like to acknowledge the fellowship and amenities provided by the Council of Scientific and Industrial Research, India (CSIR) and BITS-Pilani, Hyderabad. The second author would like to thank IIT-Guwahati for providing the fellowship and amenities. The third author was supported by BITS-Pilani, Hyderabad.


\begin{thebibliography}{99}

\bibitem{Buchholz} Buchholz RH, Rathbun RL. An infinite set of Heron triangles with two rational medians. The American mathematical monthly. 1997 Feb 1;104(2):107-15.

\bibitem{Buchholz2} Buchholz RH, Stingley RP. Heron triangles with three rational medians. Rocky Mountain Journal of Mathematics. 2019 Apr;49(2):405-17.

\bibitem{Chahal} Chahal JS. Congruent numbers and elliptic curves. The American Mathematical Monthly. 2006 Apr 1;113(4):308-17.

\bibitem{Coates} Coates JH. Congruent number problem. Pure and Applied Mathematics Quarterly. 2005 Jan;1(1):14-27.

\bibitem{conrad} Conrad K. The congruent number problem. The Harvard College Mathematics Review. 2008;2(2):58-74.



\bibitem{Goins} Goins EH, Maddox D. Heron triangles via elliptic curves. The Rocky Mountain Journal of Mathematics. 2006 Jan 1:1511-1526.

\bibitem{Halbeisen}  Halbeisen L, Hungerbühler N. Heron triangles and their elliptic curves. Journal of Number Theory. 2020 Aug 1;213:232-53.

\bibitem{Izadi} Izadi FA, Khoshnam F, Nabardi K. A new family of elliptic curves with positive ranks arising from the Heron triangles. arXiv preprint arXiv:1012.5835. 2010 Dec 28.


\bibitem{Johnstone} Johnstone JA, Spearman BK. Congruent number elliptic curves with rank at least three. Canadian Mathematical Bulletin. 2010 Dec;53(4):661-6.

\bibitem{Kramer} Kramer AV, Luca F. Some remarks on Heron triangles. Az Eszterházy Károly Tanárképző Főiskola tudományos közleményei (Új sorozat 27. köt.). Tanulmányok a matematikai tudományok köréből= Acta Academiae Paedagogicae Agriensis.. 2000:25-38.





\bibitem{Rusin} Rusin DJ. Rational triangles with equal area. New York J. Math. 1998;4(1):15.


\bibitem{Silverman} Silverman JH. The arithmetic of elliptic curves. Springer Science and Business Media; 2009 Apr 20.


\bibitem{Yan} Yan XH. The Diophantine equation $(m^{2}+ n^{2})^{x}+(2mn)^{y}=(m+ n)^{2z}$. International Journal of Number Theory. 2020 Sep 23;16(08):1701-1708.



\end{thebibliography}
\end{document}